\documentclass[11pt]{article}
\usepackage{mathrsfs}
\usepackage{tipa}
\usepackage{amsfonts}
\usepackage{amssymb}
\usepackage{amsmath}
\usepackage{amsbsy}
\def\proclaim#1{\par \bigskip\noindent {\bf #1}\bgroup\it\ }
\def\endproclaim{\egroup\par\bigskip}
\allowdisplaybreaks[4]
\newbox\TempBox \newbox\TempBoxA

\def\exp {\textsf{exp}}

\def\a.s{\textsf{a. s}}

\def\and{\textsf{and}}
\def\text#1{\mbox{\rm #1}}

\def\underwiggle 1{
\ifmmode\setbox\TempBox=\hbox{$ 1$}\else\setbox\TempBox=\hbox{ 1}\fi
\setbox\TempBoxA=\hbox to \wd\TempBox{\hss\char'176\hss}
\rlap{\copy\TempBox}\smash{\lower9pt\hbox{\copy\TempBoxA}} }

\begin{document}
\title{\bf  Local linear estimator for stochastic differential equations driven by $\alpha$-stable L\'{e}vy motions }
\author {Song Yu-Ping, Lin Zheng-Yan \footnote{Corresponding author, zlin@zju.edu.cn}\\
         Department of Mathematics, Zhejiang University, \\Hangzhou,
         China, 310027}
\maketitle

\date{}
 $\mathbf{Abstract}$: We study the local linear estimator for the drift coefficient of
stochastic differential equations driven by $\alpha$-stable L\'{e}vy
motions observed at discrete instants letting $T \rightarrow
\infty$. Under regular conditions, we derive the weak consistency
and central limit theorem of the estimator. Compare with
Nadaraya-Watson estimator, the local linear estimator has a bias
reduction whether kernel function is symmetric or not under
different schemes.

$\mathbf{Keyword}$: local linear estimator; stable L\'{e}vy motions;
bias reduction; consistency; central limit theorem.

$\mathbf{Mathematics~Subject~Classification}$: 60J52; 62G20; 62M05;
65C30.
\section{\bf{Introduction}}
 Continuous-time models play an important role in the study of
financial time series. Especially, many models in economics and
finance, like those for an interest rate or an asset price involve
continuous-time diffusion processes. Particularly, their theoretical
and empirical applications to finance  are quite extensive (see
Jacod and Shiryaev [18]). However, growing evidence shows that
stochastic processes with jumps are becoming more and more important
(see Andersen et al. [2]; Baskshi et al. [7]; Duffie et al. [11]).
Recently, stochastic processes with jumps as an intension of
continuous-path ones have been studied by more and more
statisticians since the financial phenomena can be better
characterized (see A$\ddot{1}$t-Sahalia and Jacod [1]; Bandi and
Nguyen [3]).

A diffusion model with continuous paths is represented by the
following stochastic differential equation:
$$dX_{t}=\mu(X_{t})dt+\sigma(X_{t})dW_{t},\eqno(1.1)$$
where $W_{t}$ is a standard Brownian motion, $\mu: \mathbb{R}
\rightarrow \mathbb{R}$ is an unknown measurable function and
$\sigma: \mathbb{R} \rightarrow \mathbb{R}_{+}$ is an unknown
positive function. Many authors have investigated nonparametric
estimations for the drift function $\mu(x)$ and the diffusion
function $\sigma(x)$, which to some extend prevent the
misspecification of the model (1.1) compare with the parametric
estimations. Prakasa Rao [27] constructed a non-parametric estimator
similar as the Nadaraya-Watson estimator for $\mu(x)$. Bandi and
Phillips [4] discuss the Nadaraya-Watson estimator for these
functions of non-stationary recurrent diffusion processes. Fan and
Zhang [15] proposed local linear estimators for them and obtained
bias reduction properties. In a finite sample, Xu [34]extended
re-weighted idea proposed by Hall and Presnell [16] to estimate
$\sigma(x)$ under recurrence. Xu [33] discussed the empirical
likelihood-based inference for nonparametric recurrent diffusions to
construct confidence intervals. Furthermore, Bandi and Phillips [5]
proposed a simple and robust approach to specify a parameter class
of diffusions and estimate the parameters of interest by minimizing
criteria based on the integrated squared difference between kernel
estimates of the drift and diffusion functions and their parametric
counterparts.

Recently, stochastic processes with jumps have been paid more
attention in various applications, for instance, financial time
series to reflect discontinuity of asset return (see Baskshi et al.
[7]; Duffie et al. [11]; Johannes [20]; Bandi and Nguyen [5]). In
this paper, we consider the stochastic process with jumps through
the stochastic differential equation driven by an $\alpha$-stable
L\'{e}vy motion (1 $< \alpha < 2$):
$$d X_{t} = \mu(X_{t-}) dt + \sigma(X_{t-}) d Z_{t},~~~~~X_{0} = \eta , \eqno(1.2)$$
where $\{Z_{t}, t \geq 0\}$ is a standard $\alpha$-stable L\'{e}vy
motion defined on a probability space $(\Omega, \mathscr{F}, P)$
equipped with a right continuous and increasing family of
$\sigma$-algebras $\{\mathscr{F}_{t}, t \geq 0\}$ and $\eta$ is a
random variable independent of $\{Z_{t}\}$. $Z_{1}$ has a $\alpha$-
stable distribution $S_{\alpha}(1, \beta, 0)$ with the
characteristic function: $$ E \exp \{i u Z_{1}\} = \exp
\left\{-|u|^{\alpha}\left(1 - i \beta sgn(u) \tan
\frac{\alpha\pi}{2}\right)\right\}, u \in \mathbb{R},\eqno(1.3)$$
\noindent where $\beta \in [-1, 1]$ is the skewness parameter. One
can refer to Sato [30], Barndorff-Nielsen et al. [6] for more
detailed properties on stable distributions. Usually, we get
observations $\{X_{t_{i}}, t_{i} = i \Delta_{n}, i = 0, 1, \ldots,
n\}$ for model (1.2), where $\Delta_{n}$ is the time frequency for
observation and $n$ is the sample size. This paper is devoted to the
nonparametric estimation of the unknown drift function. Our
estimation procedure for model (1.2) should be based on
$\{X_{t_{i}}, t_{i} = i \Delta_{n}, i = 0, 1, \ldots, n\}.$

The stochastic differential equation driven by L\'{e}vy motion has
received growing interest from both theoreticians and practitioners
recently, such as applications to finance, climate dynamics et al..
Masuda ([23], [24]) proved some probabilistic properties of a
multidimensional diffusion processes with jumps and provided mild
regularity conditions for a multidimensional Ornstein-Uhlenbeck
process driven by a general L\'{e}vy process for any initial
distribution to be exponential $\beta$- mixing. When model (1.2) is
specially a mean-reverted Ornstein-Uhlenbeck process driven by a
L\'{e}vy process, i.e. $\mu(x)$ is known to be linear with the form
$\mu(x) = \gamma - \lambda x$ and $\sigma = 1,$ where $(\gamma,
\lambda)$ is unknown parameters to be estimated. Based on
$\{X_{t_{i}}, t_{i} = i \Delta_{n}, i = 0, 1, \ldots, n\},$ Hu and
Long [17] studied the least-squares estimator for $\lambda > 0$ when
Z is symmetric $\alpha$-stable and $\gamma = 0$. Masuda
[25]considered an approximate self-weighted least absolute deviation
type estimator for $(\gamma, \lambda).$ Zhou and Yu [37]proved the
asymptotic distributions of the least squares estimator of the mean
reversion parameter $\lambda$ allowing for nonlinearity in the
diffusion function under three sampling schemes. However, in model
(1.2), the drift function $\mu(x)$ is seldom known and the diffusion
function may be nonlinear in reality. With no prior specified form
of the drift function, Long and Qian [22] discussed the
Nadaraya-Watson estimator for it and obtained the weak consistency
and central limit theorem.

The Nadaraya-Watson estimator given for $\mu(x)$ is locally
approximating $\mu(x)$ by a constant (a zero-degree polynomial).
However, in the context of nonparametric estimator with
finite-dimensional auxiliary variables, local polynomial smoothing
has become the ¡°golden standard¡± (see Fan [12], Wand and Jones
[36]). The local polynomial estimator is known to share the
simplicity and consistency of the kernel estimators as
Nadaraya-Watson or Gasser-M\"{u}ller estimators but avoids boundary
effects, at least when convergence rates are concerned. Local
polynomial smoothing at a point x fits a polynomial to the pairs
$(X_{i}; Y_{i})$ for those $X_{i}$ falling in a neighborhood of $x$
determined by a smoothing parameter $h$. The local polynomial
estimator has received increasing attention and it has gained
acceptance as an attractive method of nonparametric estimation
function and its derivatives. This smoothing method has become a
powerful and useful diagnostic tool for data analysis. In
particular, the local linear estimator locally fits a polynomial of
degree one. In this paper, we propose the local linear estimators
for drift function in model (1.2). As a nonparametric methodology,
local polynomial estimator makes use of the observation information
to estimate corresponding functions not assuming the function form.
The estimator is obtained by locally fitting a polynomial of degree
one to the data via weighted least squares and it shows advantages
compared with Nadaraya-Watson approach (see Fan and Gijbels [13]).
For further motivation and study of the local linear estimator, see
Fan and Gijbels [14], Ruppert and Wand [29], Stone [32], Cleveland
[10].

The remainder of this paper is organized as follows. In Section 2,
local linear estimator and appropriate assumptions for model (1.2)
are introduced. In Section 3, we present some technical lemmas and
asymptotic results. The proofs will be collected in Section 4.

\section{Local Linear Estimator and Assumptions}

\noindent We lay out some notations. For simplify, $X_{i}$ denotes
$X_{t_{i}}$ and we shall omit the subscript $n$ in the notation if
no confusion will be caused. We will use notation
``$\stackrel{p}{\rightarrow}$'' to denote ``convergence in
probability'', notation ``$\stackrel{a.s.}{\rightarrow}$'' to denote
``convergence almost surely'' and notation ``$\Rightarrow$'' to
denote ``convergence in distribution''.
\medskip

\noindent Local polynomial estimator firstly introduced in Fan [12]
has been widely used in regression analysis and time series
analysis. It has gained acceptance as an attractive method of
nonparametric estimation of regression function and its derivatives.
The estimator is obtained by locally fitting $p$-th polynomial to
the data via weighted least squares and it shows advantages compared
with other kernel nonparametric regression estimators. The idea of
weighted local polynomial regression is the following: under some
smoothness conditions of the curve $m(x)$, we can expand $m(x)$ in a
neighborhood of the point $x_{0}$ as follows:
\begin{eqnarray*}
m(x) & \approx & m(x_{0}) + m^{'}(x_{0})(x - x_{0}) +
\frac{m^{''}(x_{0})}{2!}(x - x_{0})^{2} + \cdots +
\frac{m^{(p)}(x_{0})}{p~!}(x - x_{0})^{p}\\
& \equiv & \sum_{j=0}^{p}\beta_{j}(x - x_{0})^{j},
\end{eqnarray*}
where $\beta_{j} = \frac{m^{(j)}(x_{0})}{j~!}.$

Thus, the problem of estimating infinite dimensional $m(x)$ is
equivalent to estimating the $p$-dimensional parameter $\beta_{0},
\beta_{0}, \cdots, \beta_{p}.$ Consider a weighted local polynomial
regression:
$$\arg \min_{\beta_{0},\beta_{1},\cdots,\beta_{p}} \sum_{i=0}^{n-1}
\Big\{Y_{i} - \sum_{j=0}^{p}\beta_{j}(X_{i} -
x)^{j}\Big\}^{2}K_{h_{n}}(X_{i} - x),$$ where $Y_{i} = \frac{X_{i+1}
- X_{i}}{\Delta}$ and
$K_{h_{n}}(\cdot)=\frac{1}{h_{n}}K(\frac{\cdot}{h_{n}}).$ is kernel
function with $h_{n}$ the bandwidth.

What we are interested in is to estimate $\hat{\mu}(x) =
\hat{\beta}_{0}$, hence as Fan and Gijbels [14] remarked, it is
reasonable for us to discuss $p = 1$: the local linear estimator for
the drift function $\mu(\cdot)$ in this paper. The local linear
estimator for $\mu(x)$ is the solution $\beta_{0}$ of the optimal
problem:

$$\arg \min_{\beta_{0},\beta_{1}} \sum_{i=0}^{n-1}
\Big\{Y_{i} - \sum_{j=0}^{1}\beta_{j}(X_{i} -
x)^{j}\Big\}^{2}K_{h}(X_{i} - x).$$

\noindent The solution of $\beta_{0}$ is
$$\hat{\mu}(x) = \frac{\sum\limits_{i=0}^{n-1}K_{h}(X_{i} - x)
\left\{\frac{S_{n,2}}{h^{2}} - (\frac{X_{i} -
x}{h})\frac{S_{n1}}{h}\right\}(X_{i+1} -
X_{i})}{\Delta\sum\limits_{i=0}^{n-1}K_{h}(X_{i} -
x)\left\{\frac{S_{n,2}}{h^{2}} - (\frac{X_{i} -
x}{h})\frac{S_{n1}}{h}\right\}},$$ where $S_{n,k} =
\sum\limits_{i=0}^{n-1}{K}_{h}(X_{i} - x)(X_{i} - x)^{k} , k=1,2.$

We can also write $\hat{\mu}(x) =
\frac{\hat{g}_{n}(x)}{\hat{h}_{n}(x)},$
$$\hat{h}_{n}(x) =
\frac{1}{n}\sum\limits_{i=0}^{n-1}K_{h}(X_{i} - x)
\left\{\frac{S_{n,2}}{nh^{2}} - \left(\frac{X_{i} -
x}{h}\right)\frac{S_{n1}}{nh}\right\},$$
$$\hat{g}_{n}(x) =
\frac{1}{n\Delta}\sum\limits_{i=0}^{n-1}K_{h}(X_{i} - x)
\left\{\frac{S_{n,2}}{nh^{2}} - \left(\frac{X_{i} -
x}{h}\right)\frac{S_{n1}}{nh}\right\}(X_{i+1} - X_{i}).$$

\bigskip

\noindent We now present some assumptions used in this paper.
\medskip

$\bf{(A.1).}$ The drift function $\mu(\cdot)$ is twice continuously
differentiable with bounded first and second order derivatives; the
diffusion function $\sigma(\cdot)$ satisfy a global Lipschitz
condition, there exists a positive constant C $>$ 0 such that
$$|\sigma(y) - \sigma(x)| \leq C |y - x|,~ y,~x~\in \mathbb{R}.$$

$\bf{(A.2).}$ There exist positive constants $\sigma_{0}$ and
$\sigma_{1}$ such that $0 < \sigma_{0} \leq \sigma(x) \leq
\sigma_{1}$ for each $x~\in \mathbb{R}.$

$\bf{(A.3).}$ The solution $X_{t}$ admits a unique invariant
distribution $F(x)$ and is geometrically strong mixing, i.e. there
exists $c_{0} > 0$ and $\rho \in (0 , 1)$ such that $\alpha_{X}(t)
\leq c_{0}\rho^{t},~t \geq 0.$

$\bf{(A.4).}$ The density function $f(x)$ of the stationary
distribution $F(x)$ is continuously differentiable and $f(x)$ $>$ 0.

$\bf{(A.5).}$ The kernel function $K(\cdot)$ is nonnegative
probability density function with compact support satisfying: $K_{2}
:= \int_{-\infty}^{+\infty}u^{2}K(u)du <
\infty,~~~\int_{-\infty}^{+\infty}K^{2}(u)du < \infty.$

$\bf{(A.6).}$ As $n \rightarrow \infty,~ h \rightarrow 0,~ \Delta
\rightarrow 0,~ \and~ nh\Delta \rightarrow \infty.$
\medskip

$\bf{Remark 2.1.}$ The condition (A.1) ensures that (1.2) admits a
unique non-plosive c\`{a}dl\`{a}g adapted solution, see Jacod and
Shiryaev [18]. (A.3) implies $X_{t}$ is ergodic and stationary. The
mixing property of a stochastic process describes the temporal
dependence in data. One can refer to Bradley [9] for different kinds
of mixing properties. For some sufficient conditions which guarantee
(A.3), one can refer to Masuda [24]. The kernel function is not
necessarily to be symmetric. Sometimes, unilateral kernel function
may make predictor easier (see Fan and Zhang [15]).
\medskip

\section{Some Technical Lemmas and Asymptotic Results}

We say that a continuous function $G~:[0, \infty) \rightarrow [0,
\infty)$ grows more slowly than $u^{\alpha}$ ($\alpha > 0$) if there
exist positive constants $c, \lambda_{0} \and \alpha_{0} < \alpha$
such that $G(\lambda u) \leq c\lambda^{\alpha_{0}}G(u)$ for all $u
> 0$ and all $\lambda \geq \lambda_{0}$.

\noindent $\bf{Lemma~ 3.1.}$ Let $\phi(t)$ be a predictable process
satisfying $\int_{0}^{T}|\phi(t)|^{\alpha}dt < \infty$ almost surely
for $T < \infty.$ We assume that either $\phi$ is nonnegative or $Z$
is symmetric. If $G(u)$ grows more slowly than $u^{\alpha}$, then
there exist positive constants $c_{1}$ and $c_{2}$ depending only on
$\alpha, \alpha_{0}, c ~and \lambda_{0}$ such that for each $T
> 0$
$$c_{1}E[G((\int_{0}^{T}|\phi(t)|^{\alpha}dt)^{1/{\alpha}})] \leq
 E[G(\sup_{t \leq T}|\int_{0}^{t}\phi(s)dZ_{s}|)]
 \leq c_{2}E[G((\int_{0}^{T}|\phi(t)|^{\alpha}dt)^{1/{\alpha}})].$$

\noindent $\bf{Remark~ 3.1.}$ This lemma can be viewed as a
generalization of Theorem 3.2 in Rosonski and Woyczynski [28], where
they only dealt with the case that $Z$ is symmetric.
\medskip

\noindent $\bf{Lemma~ 3.2.}$ Suppose that there is a deterministic
and nonnegative function $\Phi$ such that
$$\Phi^{\alpha}(T)\int_{0}^{T}|\phi(t)|^{\alpha}dt \stackrel{p}{\rightarrow}
1 ~~~as~~~ T \rightarrow \infty.$$ Then, we have
$$\Phi(T)\int_{0}^{T}|\phi(t)|dZ_{t} \Longrightarrow S_{\alpha}(1, \beta, 0).$$

\noindent $\bf{Remark~3.2.}$ This lemma can be regarded as an
extending to $\alpha$-stable case of Theorem 1.19 in Kutoyants [21].
\medskip

\noindent $\bf{Lemma~ 3.3.}$ Assumptions (A.1) - (A.6) lead to the
following result:
$$\frac{1}{n}\sum_{i=0}^{n-1}K_{h}(X_{i} - x)\left(\frac{X_{i} - x}{h}\right)^{k}
\stackrel{a.s.}{\longrightarrow} f(x)\int_{- \infty}^{+
\infty}u^{k}K(u)du
$$

\noindent $\bf{Remark~3.3.}$ In Long and Qian [22], they proved a
weaker case: $\frac{1}{n}\sum\limits_{i=0}^{n-1}K_{h}(X_{i} - x)
\stackrel{p}{\longrightarrow} f(x).$ One can easily obtain
$\hat{h}_{n}(x) \stackrel{a.s.}{\longrightarrow} K_{2}f^{2}(x) -
(K_{1}f(x))^{2}$ based on this lemma.
\medskip

\noindent $\bf{Theorem~ 3.1.}$Assume that (A.1)-(A.6) hold and
$\alpha \in (1, 2)$, then $\hat{\mu}(x)
\stackrel{p}{\longrightarrow} \mu(x)$ as $n \rightarrow \infty.$
\medskip

\noindent $\bf{Theorem~ 3.2.}$ Let $\alpha \in [1, 2)$ and assume
that (A.1) - (A.6) are satisfied.

$(i)$~ If $(n \Delta h)^{1 - \frac{1}{\alpha}} h = O(1)$ and $(n
\Delta h)^{1 - \frac{1}{\alpha}} \Delta^{\frac{1}{\kappa}} = O(1)$
for some $\kappa > \alpha$, then
$$(n \Delta h)^{1 - \frac{1}{\alpha}} \Lambda(x) (\hat{\mu}(x) - \mu(x)) \Rightarrow S_{\alpha}(1, \beta, 0)$$

$(ii)$ If $(n \Delta h)^{1 - \frac{1}{\alpha}} h^{2} = O(1)$ and $(n
\Delta h)^{1 - \frac{1}{\alpha}} \Delta^{\frac{1}{\kappa}} = O(1)$
for some $\kappa > \alpha$, then
$$(n \Delta h)^{1 - \frac{1}{\alpha}} \Lambda(x) (\hat{\mu}(x) - \mu(x) - h^{2}\Gamma_{\mu}(x)) \Rightarrow S_{\alpha}(1, \beta, 0)$$
where $\Lambda(x) = \frac{[K_{2} - (K_{1})^{2}]f(x)^{1 -
\frac{1}{\alpha}}}{\sigma(x)\big(\int_{- \infty}^{+
\infty}K^{\alpha}(u)\{K_{2} - u
K_{1}\}^{\alpha}du\big)^{\frac{1}{\alpha}}},$ and $\Gamma_{\mu}(x) =
\frac{\mu^{''}(x)[(K_{2})^{2} - K_{1}K_{3}]}{2 \big(K_{2} -
(K_{1})^{2}\big)}.$
\medskip

\noindent $\bf{Remark~3.4.}$ In Long and Qian [22], they showed the
following results under the Assumptions in this paper: $(i)$~ If $(n
\Delta h)^{1 - \frac{1}{\alpha}} h = O(1)$ and $(n \Delta h)^{1 -
\frac{1}{\alpha}} \Delta^{\frac{1}{\kappa}} = O(1)$ for some $\kappa
> \alpha$, then
$$(n \Delta h)^{1 - \frac{1}{\alpha}} \Lambda(x) (\hat{\mu}(x) - \mu(x)- hK_{1}) \Rightarrow S_{\alpha}(1, \beta, 0)$$

\noindent $(ii)$ If $(n \Delta h)^{1 - \frac{1}{\alpha}} h^{2} =
O(1)$ and $(n \Delta h)^{1 - \frac{1}{\alpha}}
\Delta^{\frac{1}{\kappa}} = O(1)$ for some $\kappa > \alpha$,
$K(\cdot)$ is symmetric,, then
$$(n \Delta h)^{1 - \frac{1}{\alpha}} \Lambda(x) (\hat{\mu}(x) - \mu(x) - h^{2}\Gamma_{\mu}(x)) \Rightarrow S_{\alpha}(1, \beta, 0)$$
where $\Lambda(x) = \frac{f(x)^{1 -
\frac{1}{\alpha}}}{\sigma(x)\big(\int_{- \infty}^{+
\infty}K^{\alpha}(u)\big)^{\frac{1}{\alpha}}},$ and $\Gamma_{\mu}(x)
= \big[\mu^{'}(x)\frac{f^{'}(x)}{f(x)} +
\frac{1}{2}\mu^{''}(x)\big]K_{2}.$ We can easily observe that the
bias in the local linear case is smaller than the one in the
Nadaraya-Watson case in comparison to the results between this paper
and Long \& Qian [22] whether $K(\cdot)$ is symmetric or not.
Furthermore, when $\alpha = 1$, $\hat{\mu}(x)$ is inconsistent with
$\mu(x)$ easily obtained from Theorem 3.2.

\section{Proofs}

\noindent $\bf{Proof ~of ~Lemma ~3.1.}$~ See Long and Qian (Lemma
2.7).
\medskip

\noindent $\bf{Proof ~of ~Lemma ~3.2.}$~ See Long and Qian (Lemma
2.6).
\medskip

\noindent $\bf{Proof ~of ~Lemma ~3.3.}$

We first note that
\begin{eqnarray*}
& ~ &\frac{1}{n}\sum_{i=0}^{n-1}K_{h}(X_{i} - x)\left(\frac{X_{i} -
x}{h}\right)^{k} - f(x)\int_{- \infty}^{+ \infty}u^{k}K(u)du\\
& = & \frac{1}{n}\sum_{i=0}^{n-1}K_{h}(X_{i} - x)\left(\frac{X_{i} -
x}{h}\right)^{k} - \frac{1}{n}\sum_{i=0}^{n-1}E\left[K_{h}(X_{i} -
x)\left(\frac{X_{i} -
x}{h}\right)^{k}\right]\\
& ~ & + \frac{1}{n}\sum_{i=0}^{n-1}E\left[K_{h}(X_{i} -
x)\left(\frac{X_{i} - x}{h}\right)^{k}\right] - f(x)\int_{-
\infty}^{+ \infty}u^{k}K(u)du. \hskip60mm(4.1)
\end{eqnarray*}

From the stationarity of $X_{t}$, we have:
\begin{eqnarray*}
\frac{1}{n}\sum_{i=0}^{n-1}E\left[K_{h}(X_{i} - x)\left(\frac{X_{i}
- x}{h}\right)^{k}\right] & = & E\left[K_{h}(X_{1} - x)(\frac{X_{1}
- x}{h})^{k}\right]\\
& = & \int_{- \infty}^{+ \infty}K_{h}(y - x)\left(\frac{y -
x}{h}\right)^{k}f(y)dy\\
& = & \int_{- \infty}^{+ \infty}K(u)u^{k}f(x + uh)du\\
& \rightarrow & f(x)\int_{- \infty}^{+ \infty}u^{k}K(u)du.
\hskip67mm(4.2)
\end{eqnarray*}

Thus, from (4.1) and (4.2) it suffices to prove that
\begin{eqnarray*}
& ~ &\frac{1}{n}\sum_{i=0}^{n-1}K_{h}(X_{i} - x)\left(\frac{X_{i} -
x}{h}\right)^{k} - \frac{1}{n}\sum_{i=0}^{n-1}E\left[K_{h}(X_{i} -
x)\left(\frac{X_{i} - x}{h}\right)^{k}\right]\\
& = & \frac{1}{n}\sum_{i=0}^{n-1}\left\{K_{h}(X_{i} -
x)\left(\frac{X_{i} - x}{h}\right)^{k} - E\left[K_{h}(X_{i} -
x)\left(\frac{X_{i} - x}{h}\right)^{k}\right]\right\}\\
& =: & \frac{1}{n}\sum_{i=0}^{n-1} \delta_{n,i}(x) \stackrel
{a.s.}{\longrightarrow} 0. \hskip129mm(4.3)
\end{eqnarray*}

Note that $\sup\limits_{0 \leq i \leq n-1} |\delta_{n,i}(x)| \leq
C_{0}h^{-1}$ a.s. for some positive constant $C_{0} < \infty$ by the
compact support of $K(\cdot).$ Applying Theorem 1.3 (2) in Bosq [8],
we have for each integer $q \in [1 , \frac{n}{2}]$ and each
$\varepsilon > 0$
$$P\left(\frac{1}{n}\left|\sum_{i=0}^{n-1}\delta_{n,i}(x)\right|
> \varepsilon\right) \leq 4 \exp \left( -
\frac{\varepsilon^{2}q}{8\nu^{2}(q)}\right) + 22\left( 1 + \frac{4
C_{0}h^{-1}}{\varepsilon}\right)^{1/2} q
\alpha_{X}([p]\Delta),\eqno(4.4)$$ where
$$\nu^{2}(q) = \frac{2}{p^{2}} s(q) +
\frac{C_{0}h^{-1}\varepsilon}{2}$$ with $p = \frac{n}{2q}$ and
\begin{eqnarray*}
s(q) & = & \max_{0 \leq j \leq 2q -1} E[([jp] + 1
-jp)\delta_{n,[jp]+1}(x) + \delta_{n,[jp]+2}(x) + \cdots \\ & ~ & +~
\delta_{n,[(j+1)p]}(x) + ((j + 1)p - [(j + 1)p])\delta_{n,[(j +
1)p]+1}(x)]^{2}.
\end{eqnarray*}

\noindent Using the H\"{o}lder inequality and stationarity of
$X_{i}$ , one can easily obtain that $s(q) = O(p^{2} h^{-1})$.

By choosing $q = [\sqrt{n\Delta} / \sqrt{h}]$ and $p = \frac{n}{2q}
= O(\sqrt{nh} / \sqrt{\Delta})$, we get
$$\frac{\varepsilon^{2}q}{8\nu^{2}(q)} = \varepsilon^{2} \cdot O(qh) = O(\varepsilon^{2}\sqrt{n \Delta h}).\eqno(4.5)$$
Moreover, we can obtain
$$22\left( 1 + \frac{4
C_{0}h^{-1}}{\varepsilon}\right)^{1/2} q \alpha_{X}([p]\Delta) \leq
C(\varepsilon) \exp (- O(\varepsilon^{2}\sqrt{n\Delta
h}))\eqno(4.6)$$ under the mixing properties of $X_{t}$ in (A.3) and
(A.6).

\noindent (4.4), (4.5) and (4.6) imply
$$P\left(\frac{1}{n}\left|\sum_{i=0}^{n-1}\delta_{n,i}(x)\right|
> \varepsilon\right) \leq C(\varepsilon)\exp (- O(\varepsilon^{2}\sqrt{n\Delta h})).$$

\noindent Therefore, $\frac{1}{n}\sum\limits_{i=0}^{n-1}K_{h}(X_{i}
- x)\left(\frac{X_{i} - x}{h}\right)^{k} -
\frac{1}{n}\sum\limits_{i=0}^{n-1}E\left[K_{h}(X_{i} -
x)\left(\frac{X_{i} - x}{h}\right)^{k}\right] \stackrel
{a.s.}{\longrightarrow} 0$ based on Borel-Cantelli lemma and (A.6).
\medskip

\noindent $\bf{Proof ~of ~Theorem ~3.1.}$

\noindent It suffices to prove that
$$\hat{g}_{n}(x) \stackrel {p}{\longrightarrow} [K_{2}f^{2}(x) - (K_{1}f(x))^{2}]\mu(x).$$

\noindent By (1.2), we first note that
\begin{eqnarray*}
\hat{g}_{n}(x)\!\!\! & = &\!\!\!
\frac{1}{n\Delta}\sum\limits_{i=0}^{n-1}K_{h}(X_{i} - x)
\left\{\frac{S_{n,2}}{nh^{2}} - \left(\frac{X_{i} -
x}{h}\right)\frac{S_{n1}}{nh}\right\}(X_{i+1} - X_{i})\\
& = & \!\!\!\frac{1}{n\Delta}\sum\limits_{i=0}^{n-1}K_{h}(X_{i} - x)
\left\{\frac{S_{n,2}}{nh^{2}} - \left(\frac{X_{i} -
x}{h}\right)\frac{S_{n1}}{nh}\right\}\left(\int_{t_{i}}^{t_{i+1}}
\mu(X_{s-}) ds~ + ~\int_{t_{i}}^{t_{i+1}}\sigma(X_{s-}) dZ_{s})\right)\\
& = & \!\!\!\frac{1}{n\Delta}\sum\limits_{i=0}^{n-1}K_{h}(X_{i} - x)
\left\{\frac{S_{n,2}}{nh^{2}} - \left(\frac{X_{i} -
x}{h}\right)\frac{S_{n1}}{nh}\right\}\left(\mu(X_{i})\Delta +
\int_{t_{i}}^{t_{i+1}}
(\mu(X_{s-}) - \mu(X_{i})) ds + \int_{t_{i}}^{t_{i+1}}\sigma(X_{s-}) dZ_{s}\right)\\
& =: & g_{n,1}(x) + g_{n,2}(x) + g_{n,3}(x).
\end{eqnarray*}

\noindent To show the convergence of $\hat{g}_{n}(x)$, we should
prove the following three results:

(i)~~~$g_{n,1}(x) \stackrel{p}{\longrightarrow} [K_{2}f^{2}(x) -
(K_{1}f(x))^{2}]\mu(x),~ as ~ n \rightarrow \infty ~;$

(ii)~~$g_{n,2}(x) \stackrel{p}{\longrightarrow} 0,~ as ~ n
\rightarrow \infty ~;$

(iii)~$g_{n,3}(x) \stackrel{p}{\longrightarrow} 0,~ as ~ n
\rightarrow \infty ~;$

\medskip

\noindent $\bf{Proof~of~(i):}$

\begin{eqnarray*}
g_{n,1}(x) & = & \mu(x)\frac{1}{n}\sum\limits_{i=0}^{n-1}K_{h}(X_{i}
- x) \left\{\frac{S_{n,2}}{nh^{2}} - \left(\frac{X_{i} -
x}{h}\right)\frac{S_{n1}}{nh}\right\}\\
& ~ & + \frac{1}{n}\sum\limits_{i=0}^{n-1}K_{h}(X_{i} - x)
\left\{\frac{S_{n,2}}{nh^{2}} - \left(\frac{X_{i} -
x}{h}\right)\frac{S_{n1}}{nh}\right\}(\mu(X_{i}) - \mu(x))\\
& =: & A_{n,1}(x) + A_{n,2}(x). \hskip118mm(4.7)
\end{eqnarray*}

\noindent Using Lemma 3.3, it is obviously that $$A_{n,1}(x)
\stackrel{a.s.}{\longrightarrow} \mu(x)[K_{2}f^{2}(x) -
(K_{1}f(x))^{2}].\eqno(4.8)$$

\noindent By the Lipschitz property of $\mu(x)$ and the stationarity
of $X_{t}$, we have

$$|A_{n,2}(x)| \leq \frac{L}{n}\sum_{i=0}^{n-1}|X_{i} - x|K_{h}(X_{i}
- x)\left|\frac{S_{n,2}}{nh^{2}}\right| +
\frac{L}{n}\sum_{i=0}^{n-1}|X_{i} - x|K_{h}(X_{i} -
x)\left|\frac{X_{i} - x}{h}\right|\left|\frac{S_{n,1}}{nh}\right|,
\eqno(4.9)$$ where $L$ denotes the bound of the first derivative of
$\mu(x)$.
\medskip

\noindent The two components of the right part are dealt with in the
same way, so we only deal with the first one for convenience.

\begin{eqnarray*}
\frac{1}{n}\sum_{i=0}^{n-1}|X_{i} - x|K_{h}(X_{i} -
x)\left|\frac{S_{n,2}}{nh^{2}}\right|\!\!\!\!\!\! &=&\!\!\!\!\!\!
\frac{1}{n}\sum_{i=0}^{n-1}\left(|X_{i} - x|K_{h}(X_{i} -
x)\left|\frac{S_{n,2}}{nh^{2}}\right| - E\left[|X_{i} -
x|K_{h}(X_{i} -
x)\left|\frac{S_{n,2}}{nh^{2}}\right|\right]\right)\\
& &\!\!\!\!\!\!+ \frac{1}{n}\sum_{i=0}^{n-1}E\left[|X_{i} -
x|K_{h}(X_{i} -
x)\left|\frac{S_{n,2}}{nh^{2}}\right|\right].\hskip47mm(4.10)
\end{eqnarray*}

\noindent We find that $|X_{i} - x|K_{h}(X_{i} -
x)\left|\frac{S_{n,2}}{nh^{2}}\right| - E\left[|X_{i} -
x|K_{h}(X_{i} - x)\left|\frac{S_{n,2}}{nh^{2}}\right|\right]$ is
a.s. uniformly bounded for each i by Lemma 3.3 and the compact
support of $K(\cdot)$. Similar as the proof of (4.3), we can show
that
$$\frac{1}{n}\sum_{i=0}^{n-1}\left(|X_{i} - x|K_{h}(X_{i} -
x)\left|\frac{S_{n,2}}{nh^{2}}\right| - E\left[|X_{i} -
x|K_{h}(X_{i} - x)\left|\frac{S_{n,2}}{nh^{2}}\right|\right]\right)
\stackrel{p}{\longrightarrow} 0.\eqno(4.11)$$

\noindent As for the second part,

\begin{eqnarray*}
\lim_{h \rightarrow 0}\frac{1}{n}\sum_{i=0}^{n-1}E\left[|X_{i} -
x|K_{h}(X_{i} - x)\left|\frac{S_{n,2}}{nh^{2}}\right|\right] \!\!\!&
= &\!\!\! \lim_{h \rightarrow 0}E\left[|X_{1} - x|K_{h}(X_{1} -
x)\left|\frac{S_{n,2}}{nh^{2}}\right|\right]\\
& = & \lim_{h \rightarrow 0}K_{2}f(x)E\left[|X_{1} -
x|K_{h}(X_{1} - x)\right]\\
& = & \lim_{h \rightarrow 0}hK_{2}f(x)\int_{- \infty}^{+
\infty}|u|K(u)f(x + uh)du\\
& = & \lim_{h \rightarrow 0}hK_{2}f^{2}(x)\int_{- \infty}^{+
\infty}|u|K(u)du \rightarrow 0. \hskip33mm(4.12)
\end{eqnarray*}

\noindent It follows that $A_{n,2} \stackrel{p}{\longrightarrow} 0 $
as $n \rightarrow \infty$. Hence $g_{n,1}(x)
\stackrel{p}{\longrightarrow} [K_{2}f^{2}(x) -
(K_{1}f(x))^{2}]\mu(x)$ by (4.7)-(4.12).
\medskip

\noindent $\bf{Proof~of~(ii):}$

\noindent We first introduce a basic inequality for (1.2):
$$\sup_{t_{i} \leq t \leq t_{i+1}}|X_{t} - X_{t_{i}}| \leq e^{L\Delta}\left(|\mu(X_{i})|\Delta
+ \sup_{t_{i} \leq t \leq
t_{i+1}}\left|\int_{t_{i}}^{t}\sigma(X_{s-})dZ_{s}\right|\right),\eqno(4.13)$$
\noindent which one can refer to Long \& Qian [22], Shimizu \&
Yoshida [31] and Jacod \& Protter [19].

\begin{eqnarray*}
|g_{n,2}(x)| & \leq &
\frac{1}{n\Delta}\sum\limits_{i=0}^{n-1}K_{h}(X_{i} - x)
\left\{\left|\frac{S_{n,2}}{nh^{2}}\right| + \left|\frac{X_{i} -
x}{h}\right|\left|\frac{S_{n1}}{nh}\right|\right\}
\int_{t_{i}}^{t_{i+1}} |\mu(X_{s-}) - \mu(X_{i})| ds\\
& \leq & \frac{L}{n\Delta}\sum\limits_{i=0}^{n-1}K_{h}(X_{i} - x)
\left\{\left|\frac{S_{n,2}}{nh^{2}}\right| + \left|\frac{X_{i} -
x}{h}\right|\left|\frac{S_{n1}}{nh}\right|\right\}
\int_{t_{i}}^{t_{i+1}} |X_{s-} - X_{i}| ds\\
& \leq & \frac{L}{n}\sum\limits_{i=0}^{n-1}K_{h}(X_{i} - x)
\left\{\left|\frac{S_{n,2}}{nh^{2}}\right| + \left|\frac{X_{i} -
x}{h}\right|\left|\frac{S_{n1}}{nh}\right|\right\}\sup_{t_{i} \leq t
\leq t_{i+1}}|X_{t} - X_{t_{i}}|\\
& \leq & \frac{Le^{L\Delta}}{n}\sum\limits_{i=0}^{n-1}K_{h}(X_{i} -
x) \left\{\left|\frac{S_{n,2}}{nh^{2}}\right| + \left|\frac{X_{i} -
x}{h}\right|\left|\frac{S_{n1}}{nh}\right|\right\}\left(|\mu(X_{i})|\Delta
+ \sup_{t_{i} \leq t \leq
t_{i+1}}\left|\int_{t_{i}}^{t}\sigma(X_{s-})dZ_{s}\right|\right)\\
& =: & A_{n,3}(x) + A_{n,4}(x) + A_{n,5}(x). \hskip98mm(4.14)
\end{eqnarray*}

\noindent Similarly as the proof of (4.1), we know that
$$\frac{L\Delta e^{L\Delta}}{n}\sum\limits_{i=0}^{n-1}K_{h}(X_{i}
- x) \left\{\left|\frac{S_{n,2}}{nh^{2}}\right| + \left|\frac{X_{i}
- x}{h}\right|\left|\frac{S_{n1}}{nh}\right|\right\}|\mu(X_{i})| -
Q(x) \stackrel{p}{\rightarrow} 0,$$ where $Q(x) = L\Delta
e^{L\Delta}|\mu(x)|\left[K_{2}f^{2}(x) + \left(\int_{- \infty}^{+
\infty}|u|K(u)du\right)^{2}f^{2}(x)\right] \rightarrow 0,$

\noindent that is $$\frac{L\Delta
e^{L\Delta}}{n}\sum\limits_{i=0}^{n-1}K_{h}(X_{i} - x)
\left\{\left|\frac{S_{n,2}}{nh^{2}}\right| + \left|\frac{X_{i} -
x}{h}\right|\left|\frac{S_{n1}}{nh}\right|\right\}|\mu(X_{i})|
\stackrel{p}{\rightarrow} 0.\eqno(4.15)$$
\medskip

\noindent$A_{n,4}(x), A_{n,5}(x)$ are dealt with the same approach,
hence here we only verify $A_{n,4}(x) \stackrel{p}{\rightarrow} 0.$
\medskip

\noindent As for $A_{n,4}(x)$, according to Lemma 3.3, we only need
to verify
$$\frac{1}{n}\sum\limits_{i=0}^{n-1}K_{h}(X_{i} -
x)\sup_{t_{i} \leq t \leq
t_{i+1}}\left|\int_{t_{i}}^{t}\sigma(X_{s-})dZ_{s}\right|
\stackrel{p}{\rightarrow} 0.\eqno(4.16)$$

\noindent By Markov inequality and Lemma 3.1, we have
\begin{eqnarray*}
P\left(\frac{1}{n}\sum\limits_{i=0}^{n-1}K_{h}(X_{i} - x)\sup_{t_{i}
\leq t \leq
t_{i+1}}\left|\int_{t_{i}}^{t}\sigma(X_{s-})dZ_{s}\right| >
\varepsilon\right) & \leq &
\frac{1}{n\varepsilon}\sum_{i=0}^{n-1}E\left[\sup_{t_{i} \leq t \leq
t_{i+1}}\left|\int_{t_{i}}^{t}K_{h}(X_{i} -
x)\sigma(X_{s-})dZ_{s}\right|\right]\\
& \leq &
\frac{c_{1}}{n\varepsilon}\sum_{i=0}^{n-1}E\left[\left(\int_{t_{i}}^{t_{i+1}}K_{h}^{\alpha}(X_{i}
- x)\sigma^{\alpha}(X_{s-})ds\right)^{1/\alpha}\right]\\
& \leq &
\frac{c_{1}}{n\varepsilon}\sum_{i=0}^{n-1}E\left[K_{h}(X_{i} - x)\sigma_{1}\Delta^{1/\alpha}\right]\\
& \leq & O(\Delta^{1/\alpha}) \rightarrow 0.\hskip49mm(4.17)\\
\end{eqnarray*}

\noindent Now, $g_{n,2}(x) \stackrel{p}{\longrightarrow} 0$ by
(4.14), (4.15) and (4.16).
\medskip

\noindent $\bf{Proof~of~(iii):}$

\begin{eqnarray*}
|g_{n,3}(x)| & = &
\left|\frac{1}{n\Delta}\sum\limits_{i=0}^{n-1}K_{h}(X_{i} - x)
\left\{\frac{S_{n,2}}{nh^{2}} - \left(\frac{X_{i} -
x}{h}\right)\frac{S_{n1}}{nh}\right\}\left(
\int_{t_{i}}^{t_{i+1}}\sigma(X_{s-}) dZ_{s}\right)\right|\\
& \leq & \left|\frac{S_{n,2}}{nh^{2}}\right|\cdot
\frac{1}{n\Delta}\left|\sum\limits_{i=0}^{n-1}K_{h}(X_{i} -
x)\int_{t_{i}}^{t_{i+1}}\sigma(X_{s-}) dZ_{s}\right|\\ & ~ & +
\left|\frac{S_{n,1}}{nh}\right|\cdot
\frac{1}{n\Delta}\left|\sum\limits_{i=0}^{n-1}K_{h}(X_{i} -
x)\left(\frac{X_{i} -
x}{h}\right)\int_{t_{i}}^{t_{i+1}}\sigma(X_{s-})
dZ_{s}\right|.\hskip45mm(4.18)
\end{eqnarray*}

\noindent By Lemma 3.3, we only need to prove
$$\frac{1}{n\Delta}\left|\sum\limits_{i=0}^{n-1}K_{h}(X_{i} -
x)\int_{t_{i}}^{t_{i+1}}\sigma(X_{s-}) dZ_{s}\right|
\stackrel{p}{\rightarrow} 0 \eqno(4.19)$$ and
$$\frac{1}{n\Delta}\left|\sum\limits_{i=0}^{n-1}K_{h}(X_{i} -
x)\left(\frac{X_{i} -
x}{h}\right)\int_{t_{i}}^{t_{i+1}}\sigma(X_{s-}) dZ_{s}\right|
\stackrel{p}{\rightarrow} 0.\eqno(4.20)$$

\noindent We only prove (4.19) for simplicity. Denote
$$\phi_{n,1}(t, x) = \sum_{i=0}^{n-1}\frac{1}{h^{1/\alpha}}K\left(\frac{X_{i} -x}{h}\right)
\sigma(X_{t-}) {1}_{(t_{i}, t_{i+1}]}(t),$$ we have
$\frac{1}{n\Delta}\left|\sum\limits_{i=0}^{n-1}K_{h}(X_{i} -
x)\int_{t_{i}}^{t_{i+1}}\sigma(X_{s-}) dZ_{s}\right| =
\frac{1}{n\Delta h^{\frac{\alpha -
1}{\alpha}}}\left|\int_{0}^{t_{n}}\phi_{n,1}(t, x) dZ_{t}\right|.$

\noindent With the same argument as the proof of $A_{n, 4}(x)$, we
have that:
$$P\left(\frac{1}{n\Delta}\left|\sum\limits_{i=0}^{n-1}K_{h}(X_{i} -
x)\int_{t_{i}}^{t_{i+1}}\sigma(X_{s-}) dZ_{s}\right| >
\varepsilon\right) = O \left((n\Delta
h)^{\frac{1-\alpha}{\alpha}}\right).\eqno(4.21)$$

\noindent We get $g_{n,3}(x) \stackrel{p}{\longrightarrow} 0$ by
(4.18)-(4.21) and Assumption (A.6).
\bigskip

\noindent $\bf{Proof ~of ~Theorem ~3.2.}$

\noindent Note that
$$(n\Delta h)^{1 - \frac{1}{\alpha}}\Lambda(x) (\hat{\mu}(x) - \mu(x)) =
\frac{(n\Delta h)^{1 -
\frac{1}{\alpha}}\Lambda(x)\left[\hat{g}_{n}(x) -
\mu(x)\hat{h}_{n}(x)\right]}{\hat{h}_{n}(x)} =:
\frac{B_{n}(x)}{\hat{h}_{n}(x)}.$$

\noindent We have obtained $\hat{h}_{n}(x) \rightarrow [K_{2} -
K_{1}^{2}]f^{2}(x)$ applying Lemma 3.3, so it is enough to study the
asymptotic behavior of $B_{n}(x)$.

\begin{eqnarray*}
B_{n}(x) & = & (n\Delta h)^{1 -
\frac{1}{\alpha}}\Lambda(x)\left[g_{n,1}(x) -
\mu(x)\hat{h}_{n}(x)\right] + (n\Delta h)^{1 -
\frac{1}{\alpha}}\Lambda(x)g_{n,2}(x) + (n\Delta h)^{1 -
\frac{1}{\alpha}}\Lambda(x)g_{n,3}(x)\\
& =: & B_{n,1}(x) + B_{n,2}(x) + B_{n,3}(x).\hskip99mm(4.22)
\end{eqnarray*}

\noindent $\bf{Proof~of~ B_{n,1}(x):}$

\noindent We can express $ B_{n,1}(x)$ as
$$ B_{n,1}(x) = (n\Delta h)^{1 -
\frac{1}{\alpha}}\Lambda(x)
\frac{1}{n}\sum\limits_{i=0}^{n-1}K_{h}(X_{i} - x)
\left\{\frac{S_{n,2}}{nh^{2}} - \left(\frac{X_{i} -
x}{h}\right)\frac{S_{n1}}{nh}\right\} (\mu(X_{i}) - \mu(x)).$$

\noindent Using Taylor's expansion, we get
$$\mu(X_{i}) - \mu(x) = \mu^{'}(x)(X_{i} - x) + \frac{1}{2}\mu^{''}(x + \theta_{i}(X_{i} - x))(X_{i} - x)^{2},$$
where $\theta_{i}$ is some random variable satisfying $\theta_{i}
\in [0, 1].$

\noindent Under a simple calculus, we have
$\sum\limits_{i=0}^{n-1}K_{h}(X_{i} - x)
\left\{\frac{S_{n,2}}{nh^{2}} - \left(\frac{X_{i} -
x}{h}\right)\frac{S_{n1}}{nh}\right\}\cdot (X_{i}- x) \equiv 0,$ so

\begin{eqnarray*}
B_{n,1}(x) & = & (n\Delta h)^{1 - \frac{1}{\alpha}}\Lambda(x)
\frac{1}{n}\sum\limits_{i=0}^{n-1}K_{h}(X_{i} - x)
\left\{\frac{S_{n,2}}{nh^{2}} - \left(\frac{X_{i} -
x}{h}\right)\frac{S_{n1}}{nh}\right\} \frac{1}{2}\mu^{''}(x)(X_{i} -
x)^{2}\\
& ~ & + (n\Delta h)^{1 - \frac{1}{\alpha}}\Lambda(x)
\frac{1}{n}\sum\limits_{i=0}^{n-1}K_{h}(X_{i} - x)
\left\{\frac{S_{n,2}}{nh^{2}} - \left(\frac{X_{i} -
x}{h}\right)\frac{S_{n1}}{nh}\right\} (X_{i} -
x)^{2}\frac{1}{2}[\mu^{''}(x + \theta_{i}(X_{i} - x))
-\mu^{''}(x)]\\
& =: & B_{n,1}^{(1)}(x) + B_{n,1}^{(2)}(x). \hskip116mm(4.23)
\end{eqnarray*}

\noindent By the stationary of $X_{t}$, we get
\begin{eqnarray*}
B_{n,1}^{(1)}(x) & = & \frac{1}{2}\mu^{''}(x)(n\Delta h)^{1 -
\frac{1}{\alpha}}\Lambda(x)E\left[K_{h}(X_{1} - x)
\left\{\frac{S_{n,2}}{nh^{2}} - \left(\frac{X_{1} -
x}{h}\right)\frac{S_{n1}}{nh}\right\} (X_{1} - x)^{2}\right]\\
& ~ & + \frac{1}{2}\mu^{''}(x)(n\Delta h)^{1 -
\frac{1}{\alpha}}\Lambda(x)\frac{1}{n}\sum_{i=0}^{n-1}\Big\{K_{h}(X_{i}
- x) \left\{\frac{S_{n,2}}{nh^{2}} - \left(\frac{X_{i} -
x}{h}\right)\frac{S_{n1}}{nh}\right\} (X_{i} - x)^{2}\\ & ~ & -
E\left[K_{h}(X_{i} - x) \left\{\frac{S_{n,2}}{nh^{2}} -
\left(\frac{X_{i} - x}{h}\right)\frac{S_{n1}}{nh}\right\} (X_{i} -
x)^{2}\right]\Big\}\\
& = & D_{n,1}(x) + D_{n,2}(x). \hskip116mm(4.24)
\end{eqnarray*}

\noindent One can easily obtain
$$D_{n,1}(x) = \frac{1}{2}\mu^{''}(x)\Lambda(x)[(K_{2})^{2} - K_{1}K_{3}](f(x))^{2}(n\Delta h)^{1 -
\frac{1}{\alpha}}h^{2}(1 + o(1)).\eqno(4.25)$$

\noindent Denote
$$D_{n,2}(x) := \frac{1}{2}\mu^{''}(x)\Lambda(x)\frac{1}{n}\sum_{i=0}^{n-1}\xi_{n,i}.\eqno(4.26)$$

\noindent Note that $\sup\limits_{0 \leq i \leq n-1}|\xi_{n,i}| \leq
M_{1}(n\Delta h)^{1 - \frac{1}{\alpha}}h$ a.s. for some positive
constant $M_{1} < \infty.$

It follows from the proof of (4.3) that
$\frac{1}{n}\sum\limits_{i=0}^{n-1}\xi_{n,i}
\stackrel{p}{\rightarrow} 0$ in exponential rate. However, in
addition, we should calculate $\tilde{s}(q)$, which may be a little
different from $s(q).$

\begin{eqnarray*}
\tilde{s}(q) & = & \max_{0 \leq j \leq 2q -1} E[([jp] + 1
-jp)\xi_{n,[jp]+1}(x) + \xi_{n,[jp]+2}(x) + \cdots \\ & ~ & +~
\xi_{n,[(j+1)p]}(x) + ((j + 1)p - [(j + 1)p])\xi_{n,[(j +
1)p]+1}(x)]^{2}.
\end{eqnarray*}

\noindent Using Billingsley's inequality in Bosq [8], $\tilde{s}(q)
= O(p(n\Delta h)^{2\left(1 -
\frac{1}{\alpha}\right)}\Delta^{-1}h^{2}).$

\begin{eqnarray*}
|B_{n,1}^{(2)}(x)| & \leq & \frac{1}{2}\Lambda(x)\sup\limits_{|x -
y| \leq M h} |\mu^{''}(x) - \mu^{''}(y)|(n\Delta h)^{1 -
\frac{1}{\alpha}}\Big(\frac{1}{n}\sum\limits_{i=0}^{n-1}K_{h}(X_{i}
- x) (X_{i} - x)^{2}\left|\frac{S_{n,2}}{nh^{2}}\right|\\
& & + \frac{1}{nh}\sum\limits_{i=0}^{n-1}K_{h}(X_{i} - x)|X_{i} -
x|^{3}\left|\frac{S_{n1}}{nh}\right|\Big)\\
& = & o(1)(n\Delta h)^{1 - \frac{1}{\alpha}} h^{2} \cdot
\frac{1}{n}\sum\limits_{i=0}^{n-1}K_{h}(X_{i} - x)\\
& = & o_{p}(1) (n\Delta h)^{1 - \frac{1}{\alpha}} h^{2}.
\hskip115mm(4.27)
\end{eqnarray*}

\noindent In conclusion, it follows from (4.23)-(4.27) that
$$B_{n,1}(x) = o_{p}(1) + o_{p}(1) (n\Delta h)^{1 -
\frac{1}{\alpha}} h^{2} +
\frac{1}{2}\mu^{''}(x)\Lambda(x)[(K_{2})^{2} -
K_{1}K_{3}](f(x))^{2}(n\Delta h)^{1 - \frac{1}{\alpha}}h^{2}(1 +
o(1)).\eqno(4.28)$$

\noindent $\bf{Proof~of~ B_{n,2}(x):}$ By (4.13), we have
\begin{eqnarray*}
B_{n,2}(x) & = & (n\Delta h)^{1 -
\frac{1}{\alpha}}\Lambda(x)g_{n,2}(x)\\
& \leq & (n\Delta h)^{1 - \frac{1}{\alpha}}\Lambda(x)\Big[L\Delta
e^{L\Delta}\frac{1}{n}\sum\limits_{i=0}^{n-1}K_{h}(X_{i} - x)
\left\{\left|\frac{S_{n,2}}{nh^{2}}\right| + \left|\frac{X_{i} -
x}{h}\right|\left|\frac{S_{n1}}{nh}\right|\right\}|\mu(X_{i})|\\ & ~
&+
L\Delta^{\frac{1}{\kappa}}\frac{1}{n\Delta^{\frac{1}{\kappa}}}\sum\limits_{i=0}^{n-1}K_{h}(X_{i}
- x) \left\{\left|\frac{S_{n,2}}{nh^{2}}\right| + \left|\frac{X_{i}
- x}{h}\right|\left|\frac{S_{n1}}{nh}\right|\right\}\sup_{t_{i} \leq
t \leq
t_{i+1}}\left|\int_{t_{i}}^{t}\sigma(X_{s-})dZ_{s}\right|\Big].\hskip10mm(4.29)
\end{eqnarray*}

\noindent Similar to the proof of (4.1), under given conditions and
Lemma 3.3, we can obtain
$$\frac{1}{n}\sum\limits_{i=0}^{n-1}K_{h}(X_{i} - x)
\left\{\left|\frac{S_{n,2}}{nh^{2}}\right| + \left|\frac{X_{i} -
x}{h}\right|\left|\frac{S_{n1}}{nh}\right|\right\}|\mu(X_{i})|
\stackrel{p}{\rightarrow} \left[K_{2}f^{2}(x) +
|K_{1}|f^{2}(x)\int_{- \infty}^{+ \infty}|u|K(u)du
\right]|\mu(x)|.\eqno(4.30)$$

\noindent As in the proof of (4.17) and Lemma 3.3, we can show that
$$P\left(\frac{1}{n\Delta^{\frac{1}{\kappa}}}\sum\limits_{i=0}^{n-1}K_{h}(X_{i}
- x) \left\{\left|\frac{S_{n,2}}{nh^{2}}\right| + \left|\frac{X_{i}
- x}{h}\right|\left|\frac{S_{n1}}{nh}\right|\right\}\sup_{t_{i} \leq
t \leq t_{i+1}}\left|\int_{t_{i}}^{t}\sigma(X_{s-})dZ_{s}\right| >
\varepsilon\right) \leq O(\Delta^{\frac{1}{\alpha} -
\frac{1}{\kappa}}).\eqno(4.31)$$

\noindent It follows from (4.29)-(4.31) and $\kappa > \alpha$ that
$$B_{n,2}(x) = O_{p}(1)\cdot(n\Delta h)^{1 -
\frac{1}{\alpha}}\Delta + o_{p}(1)\cdot(n\Delta h)^{1 -
\frac{1}{\alpha}}\Delta^{\frac{1}{\kappa}}.\eqno(4.32)$$

\noindent $\bf{Proof~of~ B_{n,3}(x):}$

\noindent According to (4.18),

\begin{eqnarray*}
B_{n,3}(x) & = & (n\Delta h)^{1 -
\frac{1}{\alpha}}\Lambda(x)\frac{1}{n\Delta}\sum\limits_{i=0}^{n-1}K_{h}(X_{i}
- x) \left\{\frac{S_{n,2}}{nh^{2}} - \left(\frac{X_{i} -
x}{h}\right)\frac{S_{n1}}{nh}\right\}\left(
\int_{t_{i}}^{t_{i+1}}\sigma(X_{s-}) dZ_{s}\right)\\
& =: & (n\Delta h)^{1 - \frac{1}{\alpha}}\Lambda(x) B_{n,3}^{'}(x).
\hskip107mm(4.33)
\end{eqnarray*}

\noindent Denote
$$\phi_{n,2}(t, x) = \sum_{i=0}^{n-1}\frac{1}{h^{1/\alpha}}K\left(\frac{X_{i} -x}{h}\right)
\left\{\frac{S_{n,2}}{nh^{2}} - \left(\frac{X_{i} -
x}{h}\right)\frac{S_{n1}}{nh}\right\} \sigma(X_{t-}){1}_{(t_{i},
t_{i+1}]}(t),$$

\noindent then
$$B_{n,3}^{'}(x) = \int_{0}^{t_{n}}\phi_{n,2}(t, x) dZ_{t}.$$

\noindent Let
$$\Phi_{t_{n}} = \left[t_{n} \sigma^{\alpha}(x)f(x)\int_{- \infty}^{+ \infty}
K^{\alpha}(u)\left\{K_{2}f(x) -
uK_{1}f(x)\right\}^{\alpha}du\right]^{-\frac{1}{\alpha}}.$$

\noindent Then, we have
\begin{eqnarray*}
\Phi^{\alpha}_{t_{n}} \cdot \int_{0}^{t_{n}} \phi_{n,2}^{\alpha}(t,
x) dt & = & \Phi^{\alpha}_{t_{n}} \cdot
\sum_{i=0}^{n-1}\frac{1}{h}K^{\alpha}\left(\frac{X_{i} -x}{h}\right)
\left\{\frac{S_{n,2}}{nh^{2}} - \left(\frac{X_{i} -
x}{h}\right)\frac{S_{n1}}{nh}\right\}^{\alpha}\int_{t_{i}}^{t_{i +
1}}\sigma^{\alpha}(X_{s-})ds\\
& = & \Phi^{\alpha}_{t_{n}} \cdot
\sum_{i=0}^{n-1}\frac{1}{h}K^{\alpha}\left(\frac{X_{i} -x}{h}\right)
\left\{\frac{S_{n,2}}{nh^{2}} - \left(\frac{X_{i} -
x}{h}\right)\frac{S_{n1}}{nh}\right\}^{\alpha}\sigma^{\alpha}(X_{i})\Delta\\
& ~ & + \Phi^{\alpha}_{t_{n}} \cdot
\sum_{i=0}^{n-1}\frac{1}{h}K^{\alpha}\left(\frac{X_{i} -x}{h}\right)
\left\{\frac{S_{n,2}}{nh^{2}} - \left(\frac{X_{i} -
x}{h}\right)\frac{S_{n1}}{nh}\right\}^{\alpha}\int_{t_{i}}^{t_{i +
1}}(\sigma^{\alpha}(X_{s-}) -
\sigma^{\alpha}(X_{i}))ds\\
& =: & I + J.\hskip121mm(4.34)
\end{eqnarray*}

\noindent It follows from the proof of (4.1) that
$$\frac{1}{n}\sum_{i=0}^{n-1}\frac{1}{h}K^{\alpha}\left(\frac{X_{i} -x}{h}\right)
\left\{\frac{S_{n,2}}{nh^{2}} - \left(\frac{X_{i} -
x}{h}\right)\frac{S_{n1}}{nh}\right\}^{\alpha}\sigma^{\alpha}(X_{i})
\stackrel{p}{\rightarrow} \sigma^{\alpha}(x)f(x)\int_{- \infty}^{+
\infty} K^{\alpha}(u)\left\{K_{2}f(x) -
uK_{1}f(x)\right\}^{\alpha}du.$$

\noindent Therefore, we have
\begin{eqnarray*}
I & = & \frac{1}{\sigma^{\alpha}(x)f(x)\int_{- \infty}^{+ \infty}
K^{\alpha}(u)\left\{K_{2}f(x) - uK_{1}f(x)\right\}^{\alpha}du} \cdot
\frac{1}{n}\sum_{i=0}^{n-1}\frac{1}{h}K^{\alpha}\left(\frac{X_{i}
-x}{h}\right) \left\{\frac{S_{n,2}}{nh^{2}} - \left(\frac{X_{i} -
x}{h}\right)\frac{S_{n1}}{nh}\right\}^{\alpha}\sigma^{\alpha}(X_{i})\\
& \stackrel{p}{\rightarrow} & 1.\hskip157mm(4.35)
\end{eqnarray*}

\noindent Next we deal with the second term J. By the mean-value
theorem, the Lipschitz condition (A.1) and (4.13), we have
\begin{eqnarray*}
|J| & = &
\Phi_{t_{n}}^{\alpha}\left|\sum_{i=0}^{n-1}\frac{1}{h}K^{\alpha}\left(\frac{X_{i}
-x}{h}\right) \left\{\frac{S_{n,2}}{nh^{2}} - \left(\frac{X_{i} -
x}{h}\right)\frac{S_{n1}}{nh}\right\}^{\alpha}\int_{t_{i}}^{t_{i +
1}}(\sigma^{\alpha}(X_{s-}) - \sigma^{\alpha}(X_{i}))ds\right|\\
& \leq &
\Phi_{t_{n}}^{\alpha}\sum_{i=0}^{n-1}\frac{1}{h}K^{\alpha}\left(\frac{X_{i}
-x}{h}\right) \left\{\frac{S_{n,2}}{nh^{2}} + \left|\frac{X_{i} -
x}{h}\right|\left|\frac{S_{n1}}{nh}\right|\right\}^{\alpha}\int_{t_{i}}^{t_{i
+ 1}}|\sigma^{\alpha}(X_{s-}) - \sigma^{\alpha}(X_{i})|ds\\
& \leq &
\frac{C}{n\Delta}\sum_{i=0}^{n-1}\frac{1}{h}K^{\alpha}\left(\frac{X_{i}
-x}{h}\right) \left\{\frac{S_{n,2}}{nh^{2}} + \left|\frac{X_{i} -
x}{h}\right|\left|\frac{S_{n1}}{nh}\right|\right\}^{\alpha}\int_{t_{i}}^{t_{i
+ 1}}\alpha|\zeta^{\alpha - 1}||\sigma(X_{s-}) - \sigma(X_{i})|ds\\
& \leq &
\frac{C}{n\Delta}\sum_{i=0}^{n-1}\frac{1}{h}K^{\alpha}\left(\frac{X_{i}
-x}{h}\right) \left\{\frac{S_{n,2}}{nh^{2}} + \left|\frac{X_{i} -
x}{h}\right|\left|\frac{S_{n1}}{nh}\right|\right\}^{\alpha} \cdot
\Delta \cdot
\sup\limits_{t_{i}\leq t \leq t_{i+1}}|X_{t} -X{i}|\\
& \leq & \frac{C
e^{L\Delta}}{n}\sum_{i=0}^{n-1}\frac{1}{h}K^{\alpha}\left(\frac{X_{i}
-x}{h}\right) \left\{\frac{S_{n,2}}{nh^{2}} + \left|\frac{X_{i} -
x}{h}\right|\left|\frac{S_{n1}}{nh}\right|\right\}^{\alpha}|\mu(X_{i})|\Delta\\
& ~ & + \frac{C
e^{L\Delta}}{n}\sum_{i=0}^{n-1}\frac{1}{h}K^{\alpha}\left(\frac{X_{i}
-x}{h}\right) \left\{\frac{S_{n,2}}{nh^{2}} + \left|\frac{X_{i} -
x}{h}\right|\left|\frac{S_{n1}}{nh}\right|\right\}^{\alpha}\sup\limits_{t_{i}\leq
t \leq t_{i+1}}|\int_{t_{i}}^{t}\sigma(X_{s-})dZ_{s}|\\
& =: & J_{1} + J_{2}, \hskip141mm(4.36)
\end{eqnarray*}
where $\zeta$ is some random variable satisfying $\zeta \in
[\sigma(X_{s-}) , \sigma(X_{i})]$ or $\zeta \in [\sigma(X_{i}) ,
\sigma(X_{s-})].$

\noindent Applying Lemma 3.3, and the similar proof of Lemma 3.3, we
can obtain
$$\frac{1}{n}\sum_{i=0}^{n-1}\frac{1}{h}K^{\alpha}\left(\frac{X_{i}
-x}{h}\right) \left\{\frac{S_{n,2}}{nh^{2}} + \left|\frac{X_{i} -
x}{h}\right|\left|\frac{S_{n1}}{nh}\right|\right\}^{\alpha}|\mu(X_{i})|
\stackrel{p}{\rightarrow} |\mu(x)|f(x) \int_{- \infty}^{+ \infty}
K^{\alpha}(u)\left\{K_{2}f(x) + |u||K_{1}|f(x)\right\}^{\alpha}du,$$

\noindent hence $J_{1} \stackrel{p}{\rightarrow} 0.$ As for $J_{2}
\stackrel{p}{\rightarrow} 0$, one can refer to (4.31).

\noindent From (4.33)-(4.36) and Lemma 3.2, we have
$$\Phi_{t_{n}} B^{'}_{n,3}(x) \Rightarrow S_{\alpha}(1, \beta, 0).\eqno(4.37)$$

\noindent Therefore, it follows from (4.33), (4.37) and Lemma 3.3
that
\begin{eqnarray*}
B_{n,3}(x) & = & [K_{2}f^{2}(x) - (K_{1}f(x))^{2}] \cdot
\Phi_{t_{n}} \cdot B^{'}_{n,3}(x)\\
& \Rightarrow & [K_{2}f^{2}(x) - (K_{1}f(x))^{2}]S_{\alpha}(1,
\beta, 0).
\end{eqnarray*}

\noindent By Slutsky's theorem and remark 3.3, we find
\begin{eqnarray*}
& ~ &(n \Delta h)^{1 - \frac{1}{\alpha}} \Lambda(x) (\hat{\mu}(x) -
\mu(x) - h^{2}\Gamma_{\mu}(x))\\
& = & \frac{B_{n}(x)}{\hat{h}_{n}(x)} - (n \Delta h)^{1 -
\frac{1}{\alpha}} h^{2} \Lambda(x)\Gamma_{\mu}(x)\\
& = & \frac{B_{n}(x) - (n \Delta h)^{1 - \frac{1}{\alpha}} h^{2}
\Lambda(x)\Gamma_{\mu}(x)h(x)}{\hat{h}_{n}(x)} + (n \Delta h)^{1 -
\frac{1}{\alpha}} h^{2}
\Lambda(x)\Gamma_{\mu}(x)\left(\frac{h(x)}{\hat{h}_{n}(x)} -
1\right)\\
& \Rightarrow & S_{\alpha}(1, \beta, 0).
\end{eqnarray*}

\noindent This complete the proof of (ii) in Theorem 3.2. The proof
of (i) is similar.

\section{\bf{References}}

[1] A\"{1}t-Sahalia, Y., Jacod, J. (2009) Testing for jumps in a
discretely observed process. {\em Annals of
    Statistics} $~~~~~~~$37, 184-222. \\
{[2]} Andersen, T., Benzoni, L., Lund, J. (2002) An empirical
investigation of continuous time equity return $~~~~~~$models. {\em Journal of Finance} 57, 1239-1284.\\
{[3]} Bandi, F., Nguyen, T.(2003) On the functional estimation
of jump diffusion models. {\em Journal of Econo-$~~~~~~$metrics} 116, 293-328.\\
{[4]} Bandi, F., Phillips, P. (2003) Fully nonparametric estimation
of scalar diffusion models. {\em Econometrica} $~~~~~~~$71, 241-283.\\
{[5]} Bandi, F., Phillips, P. (2007) A simple approach to the
parametric estimation of potentially nonstation-$~~~~~~$ary
diffusions. {\em Journal of Econometrics} 137,
354-395.\\
{[6]} Barndorff-Nielsen, O.E., Mikosch, T., Resnick, S. (2001)
L\'{e}vy Processes, Theory and Applications.
$~~~~~~~$Birkh\"{a}user, Boston.\\
{[7]} Baskshi, G., Cao, Z.,Chen, Z. (1997) Empirical performance
of alternative option pricing models. {\em Journal $~~~~~~~$of Finance} 52, 2003-2049\\
{[8]} Bosq, D. (1996) Nonparametric Statistics for Stochastic
Processes, Lecture Notes in Statistics, Vol. 110, $~~~~~~~$Springer-Verlag, New York.\\
{[9]} Bradley, R. (2005) Basic properties of strong mixing
conditions. A survey and some open questions. $~~~~~~~${\em Probability Surveys} 2, 107-144.\\
{[10]} Cleveland, W.S. (1979) Robust locally weighted regression and
smoothing scatter plots. {\em Journal of  the  $~~~~~~~$American Statistical Association} 74 (368), 829-836.\\
{[11]} Duffe, D.,  Pan, J. and  Singleton, K.J. (2000) Transform
analysis and asset pricing for affine jump-$~~~~~~~$diffusions. {\em Econometrica} 68, 1343-1376.\\
{[12]}Fan, J. (1992) Design-adaptive  nonparametric  regression.
{\em Journal of  the  American  Statistical Association} $~~~~~~~$ 87,  998-1004.\\
{[13]} Fan, J. and Gijbels, I. (1992) Variable bandwidth and local
linear regression smoothers. {\em Annals  of $~~~~~~~$Statistics} 20, 2008-2036.\\
{[14]} Fan, J., Gijbels, I. (1996) Local Polynomial Modeling and Its
Applications. Chapman and Hall.\\
{[15]} Fan, J. and  Zhang, C. (2003) A re-examination of diffusion
estimators with applications to financial $~~~~~~~$model validation.
{\em Journal of  the  American  Statistical Association} 98, 118-
134.\\
{[16]} Hall, P., Presnell, B. (1999) Intentionally biased bootstrap methods.
{\em Journal of the Royal Statistical $~~~~~~~$Society Series B} 61, 143-158.\\
{[17]} Hu, Y., Long, H. (1999) Least squares estimator for
Ornstein-Uhlenbeck processes driven by $\alpha$-stable
$~~~~~~~$motions. {\em Stochastic Processes and their applications}
119, 2465-2480.\\
{[18]} Jacod, J., Shiryaev, A. N. (2003) Limit Theory for Stochastic
Processes. Springer, Berlin-Heidelberg-$~~~~~~$New York-Hong Kong-London-Milan-Paris-Tokyo.\\
{[19]} Jacod, J., Protter, P. (2012) Discretization of processes.
Springer-Verlag, Heidelberg-New York.\\
{[20]} Johannes, M.S. (2004) The economic and statistical role of
jumps to interest rates. {\em Journal of Finance} $~~~~~~~$59, 227-260.\\
{[21]} Kutoyants, Y. (2004) Statistical Inference for Ergodic
Diffusion Processes. Springer-Verlag, London, $~~~~~~~$Berlin, Heidelberg.\\
{[22]} Long, H., Qian, L. Nadaraya-Watson estimator for stochastic
processes driven by stable L\'{e}vy motions. $~~~~~~~$Submitted.\\
{[23]} Masuda, H. (2004) On multi-dimensional Ornstein-Uhlenbeck
processes driven by a general L\'{e}vy pro-$~~~~~~~$cess. {\em
Bernoulli}
10, 1-24.\\
{[24]} Masuda, H. (2007) Ergodicity and exponential $\beta$-mixing
bounds for multidimensional diffusions with $~~~~~~~$jumps. {\em
Stochastic
Processes and their Applications} 117, 35-56.\\
{[25]} Masuda, H. (2010) Approximate self-weighted LAD estimation of
discretely observed ergodic Ornstein-$~~~~~~~$Uhlenbeck processes.
{\em Electronic Journal of Statistics} 4, 525-565.\\
{[26]} Moloche, G. (2001). Local Nonparametric Estimation of
Scalar Diffusions. Mimeo, MIT, Working Paper.\\
{[27]} Prakasa Rao, B.L.S (1985). Estimation of the drift for
diffusion process. {\em Statistics} 16, 263-275.\\
{[28]} Rosinski, J., Woyczynski (1985). Moment inequalities for real
and vector $p$-stable stochastic integrals. $~~~~~~~${\em
Probability in Banach Spaces V}, Lecture Notes in Math.,
1153, 369-386 Springer, Berlin.\\
{[29]} Ruppert, D., Wand, M.P. (1994). Multivariate locally weighted
least squares regression. {\em Annals of
Statis-$~~~~~~~$tics} 22, 1346-1370.\\
{[30]} Sato, K.I. (1999). L\'{e}vy Processes and Infinitely
Divisible Distributions. Cambridge University Press, $~~~~~~~$Cambridge.\\
{[31]} Shimizu, Y., Yoshida, N. (2006) Estimation of parameters for
diffusion processes with jumps from $~~~~~~~$discrete observation.
{\em Statistical Inference of Stochastic Processes} 9, 227-277.\\
{[32]} Stone, C. (1982) Optimal global rates of convergence for
nonparametric regression. {\em Annals of Statistics} $~~~~~~~$10, 1040-1053.\\
{[33]} Xu, K., (2009) Empirical likelihood based inference for
nonparametric recurrent diffusions. {\em Journal of
$~~~~~~$Econometrics}
153, 65-82.\\
{[34]} Xu, K., (2010) Re-weighted functional estimation of diffusion
models. {\em Econometric Theory} 26, 541-563.\\
{[35]} Xu, Z. (2003) Statistical inference for diffusion process.
Ph.D. thesis, East China Normal University $~~~~~~$(Chinese version).\\
{[36]} Wand, M., Jones, C. (1995) Kernel Smoothing. Chapman and Hall, London, U.K..\\
{[37]} Zhou, Q., Yu, J. (2010) Asymptotic distributions of the least
squares estimator for diffusion process. $~~~~~~~~$Singapore Management University, Working Paper.\\

\end{document}